\documentclass[12pt]{article}
\usepackage{amsmath}
\usepackage{amssymb}
\usepackage[cp1251]{inputenc}
\usepackage[russian]{babel}
\newcommand{\il}[2]{\int\limits_{#1}^{#2}}
\newcommand{\ilp}[1]{\int\limits_{#1}^{+\infty}}

\newcommand{\ph}{\phantom{a}}
\newcommand{\phh}{\phantom{aaa}}

\newcommand{\sist}[2]{\left\{
\begin{array}{l}
{#1}\\
\ph\\
{#2}
\end{array}
\right.}

\begin{document}
MSC 34C10

\vskip 15pt

\centerline{\bf Oscillation criteria for   second order  two dimensional}

 \centerline{\bf linear systems of ordinary differential equations}
\vskip 15pt

\centerline{\bf G. A. Grigorian}
\centerline{\it Institute  of Mathematics NAS of Armenia}
\centerline{\it E -mail: mathphys2@instmath.sci.am}
\vskip 20 pt

Abstract. Some properties of global solutions of scalar Riccati equation are studied. On the basis of these properties  using the Whyburn's and Leighton - Nehary's theorems some oscillatory   criteria are proved for second order linear systems of ordinary differential equations.


\vskip 12pt

\noindent
{\bf  1. Introduction}.
Let  $p(t), \phantom{a} q(t), \phantom{a} r(t),  \phantom{a} r_1(t)$   and  $r_2(t)$ be real-valued continuous functions on $[t_0,+\infty)$ and let $ \phantom{a} p(t) > 0, \phantom{a} r_1(t) r_2(t) > 0, \phantom{a} t\ge t_0$.  Set $R(t) \equiv\left(\begin{array}{c}r(t) \phantom{aa} r_1(t)\\
- r_2(t)\phantom{a} r(t)
\end{array}
\right),
$
Consider the second order two dimensional  linear system of ordinary differential equations
$$
(p(t)\Phi')' + q(t)\Phi' + R(t) \Phi =0, \phantom{aaa} t\ge t_0. \eqno (1.1)
$$
Here  $\Phi = \Phi(t)\equiv colon\{\phi_1(t),\phantom{a}\phi_2(t)\}$ is the unknown continuously differentiable vector function on $[t_0,+\infty)$.

{\bf Remark 1.1.} {\it The system (1.1) is to be interpret as the following first order linear one
$$
\sist{\phi' = \phh \frac{1}{p(t)} \Psi,}{\Psi' = - R(t) \Phi - q(t) \Psi, \ph t \ge t_0.}
$$
Note that for every "bad" (continuous but not differentiable) $p(t)$ there exist many "contr bad" \ph $\Phi(t)$ in the sense that $p(t)\Phi'(t)$ \ph is continuously differentiable. For example, \linebreak $p(t) = 1 + |\sin t|, \ph \Phi(t) = colon \Bigl\{\il{t_0}{t}\frac{\theta_1(\tau) d \tau}{1 + |\sin \tau|}, \il{t_0}{t}\frac{\theta_2(\tau) d \tau}{1 + |\sin \tau|}\Bigr\}$
where $\theta_1(t)$ and $\theta_2(t)$ are arbitrary continuously differentiable functions on $[t_0, +\infty)$. In general for arbitrary "bad" \ph  $p(t) >~ 0, \linebreak t \ge t_0$ the vector functions $\Phi_{\theta_1, \theta_2}(t) = colon \Bigl\{\il{t_0}{t}\frac{\theta_1(\tau) d \tau}{p(\tau)}, \il{t_0}{t}\frac{\theta_2(\tau) d \tau}{p(\tau)}\Bigr\}$ are "contr bad" \ph
 for $p(t)$.}

{\bf Definition 1.1}.  {\it The system (1.1) is called oscillatory if  for its every solution  \linebreak $\Phi(t)\equiv~ colon\{\phi_1(t),\phantom{a}\phi_2(t)\}$ the functions $\phi_1(t)$   and  $\phi_2(t)$  have arbitrary large zeroes.}

Study of the questions of oscillation
 of linear systems of ordinary differential equations in particular of the system
(1.1) is an important problem of qualitative theory of differential equations.
The system (1.1) appears in various problems in technics,  in particular in   the study of feathered projectile dynamics (feathered projectile oscillation)  (see [1], p. 309).

 Whyburn ([2], p. 184, [3]) studied the conditions of oscillation for the system (1.1) in the particular case when $p(t)\equiv 1, \phantom{a} q(t)\equiv 0$
and proved the following theorem (see [2], p. 184, Theorem 4.376)

{\bf Theorem 1.1 (Whyburn).} {\it Suppose that $r(t)\equiv 0$ and $r_1(t) t_2(t) > 0$ in $[t_0, +\infty)$ and
$$
\ilp{t_0} t \hskip 3pt r_1(t)d t = \ilp{t_0} t \hskip 3pt r_2(t)d t =  \infty. \eqno (1.2)
$$
Then every solution $colon \{\phi(t), \psi(t) \}$ of the system
$$
\Phi'' + R(t) \Phi = 0, \phh t \ge t_0 \eqno (1.3)
$$
with $r_1(t)[\phi(t_0)\psi'(t_0) - \psi(t_0) \phi'(t_0)] \ge 0$ is oscillatory.
If (1.2) is replaced by

$$
\ilp{t_0}  r_1(t)d t = \ilp{t_0}  r_2(t)d t =  \infty
$$
then every solution  $\Phi(t)= colon \{\phi(t), \psi(t) \}$ of (1.3) is either oscillatory or there exists a number $t_1 \ge t_0$ such that $\phi^2(t) + \psi^2(t)$ tends monotonically to zero in $[t_1, +\infty)$ as $t \to +\infty$.}

Let $p_j(t), \ph q_j(t), \ph j=1,2, \ph  q_1(t) \ne 0, \ph t \ge t_0$ be real-valued continuously differentiable  functions on $[t_0, +\infty)$. Consider the second order linear system of ordinary differential equations
$$
\sist{\phi'' = p_1(t) \phi + q_1(t) \psi,}{\psi'' = p_2(t) \phi + q_2(t) \psi, \ph t \ge t_0.} \eqno (1.4)
$$
In [4] H. Kh. Abdullah proved the following oscillation theorem

{\bf Theorem 1.2 (H. Kh. Abdullah)}. {\it If there exists a real number $K\ne 0$ such that
$$q_1(t) K^2 + (q_2(t) - p_1(t))K - p_2(t) \equiv 0 \eqno (1.5)$$ and
$$
\ilp{t_0}[K q_1(t) + q_2(t)] d t = + \infty (\mbox{or} - \infty)
$$
$$
\ilp{t_0}[K p_1(t) + p_2(t)] d t = + \infty (\mbox{or} - \infty)
$$
then the system (1.4) is oscillatory.}

For  twice continuously differentiable function  $r(t)/r_1(t)$, as is shown in [2, p. 183], the system  $(1.3)$ is reducible to the following linear differential equation of fourth order
$$
[a(t)\phi'']'' + [b(t)\phi']' + c(t) \phi = 0, \phh t \ge t_0, \eqno (1.6)
$$
where $a(t) = 1/ r(t), \ph b(t)= 2 r(t)/ r_1(t), \ph c(t) = (r(t)/ r_1(t))'' + r_2(r)/ r_1(t) + r_2(t), t \ge t_0$. In the particular case $r(t)\equiv 0$ Eq (1.6) takes the form
$$
\Bigl[\frac{1}{r_1(t)}\phi''\Bigr]'' + r_2(t) \phi = 0, \phh t \ge t_0. \eqno (1.7)
$$
For this equation Leighton and Nehary have obtained the following result (see [2], p. 121 Theorem 3.24)

{\bf Theorem 1.3 (Leighton and Nehary).} {\it If $\phi_1(t)$ and $\phi_2(t)$ are non trivial solutions of Eq. (1.7) the number of zeroes of $\phi_1(t)$ on any closed interval $[t_1, t_2]$ cannot differ by more than $4$ from the number of zeroes of $\phi_2(t)$ on $[t_1, t_2]$. In particular the solutions of (1.7) are either oscillatory (i. e. each of them have arbitrary large zeroes) or nonoscillatory (i. e. each of them  have no more than a finite number of zeroes).}

Obviously from Theorem 1.1 and Theorem 1.3  it follows immediately

{\bf Theorem 1.4 (Whiburn, Leghton and Nehary).} {\it Suppose that $r(t) \equiv 0$ and $r_1(t) r_2(t) > 0, \ph t \ge t_0$ and that (1. 2) is satisfied. Then every solution of Eq. (1.7) is oscillatory.}

It should be noted that to  study the questions of oscillation of solutions of linear differential equations of fourth order  many works  are devoted (see  [3] and cited works therein,   [5 - 10]).

In section 2  some properties of global solutions of  scalar Riccati equation are studied. By use of these properties on the basis of Theorem 1.4  in  section 3  oscillatory criteria for the system (1.1) are proved.

{\bf  2.  Auxiliary propositions}.
In what follows we will assume that all solutions of the considering equations and systems of equations are real-valued.
Let  $a(t),\phantom{a} b(t), \phantom{a} c(t), \phantom{a} a_1(t),\linebreak b_1(t), \phantom{a} c_1(t)$                    be real-valued continuous functions on $[t_0,+\infty)$. Consider the Riccati equations.
$$
y'(t) + a(t) y^2(t) + b(t) y(t) + c(t) = 0, \phantom{aaa} t\ge t_0, \eqno (2.1)
$$
$$
y'(t) + a_1(t) y^2(t) + b_1(t) y(t) + c_1(t) = 0, \phantom{aaa} t\ge t_0, \eqno (2.2)
$$
Along with these equations consider the differential inequalities
$$
\eta'(t) + a(t) \eta^2(t) + b(t) \eta(t) + c(t) \ge 0, \phantom{aaa} t\ge t_0, \eqno (2.3)
$$
$$
\eta'(t) + a_1(t) \eta^2(t) + b_1(t) \eta(t) + c_1(t) \ge 0, \phantom{aaa} t\ge t_0, \eqno (2.4)
$$
For  $a(t) \ge 0 \phantom{a} (a_1(t) \ge 0), \phantom{a} t\ge t_0,$ the inequality  (2.3)   ((2.4))  has a solution on  $[t_0,+\infty)$,  satisfying any real initial condition (see [11]).

{\bf Theorem 2.1}.   {\it Let $y_0(t)$  be a solution of Eq. (2.1) on  $[t_1,+\infty) \phantom{a} (t_1 \ge t_0)$,   $\eta_0(t)$   and $\eta_1(t)$  be solutions of the inequalities (2.3) and (2.4) respectively  with $\eta_k(t_1) \ge y_0(t_1), \linebreak  k=~0,1$, and let  $a_1(t) \ge 0,$
$$
\lambda - y_0(t) + \int\limits_{t_1}^t \exp\biggl\{\int\limits_{t_1}^\tau[a_1(\xi)(\eta_0(\xi) + \eta_1(\xi)) + b_1(\xi)]d\xi\biggr\}\times \phantom{aaaaaaaaaaaaaaaaaaaaaaaaaaaaaa}
$$
$$
\phantom{aaaaaaaaa}\times [(a(t) - a_1(t)) y_0^2(t) + (b(t) - b_1(t)) y_0(t) + c(t) - c_1(t)]d\tau \ge 0, \phantom{aaa} t\ge t_0,
$$
for some $\lambda\in [y_0(t_1), \eta_0(t_1)]$. Then for every $y_{(0)} \ge y_0(t_1)$
Eq. (2.2) has a solution $y_1(t)$ on  $[t_1,+\infty)$,
satisfying the initial condition  $y_1(t_1) = y_{(0)}$,   and  $y_1(t) \ge y_0(t), \phantom{a} t\ge t_1.$}

See the proof in [12].

{\bf Theorem 2.2}. {\it Assume  $a(t) \ge 0, \phantom{a} c(t) \le 0, \phantom{a} t\ge t_0$.
Then for every $y_{(0)} \ge 0$ Eq.~ (2.1) has a solution $y_0(t)$  on $[t_0,+\infty)$, satisfying the initial condition  $y_0(t_0) = y_{(0)}$ and $y_0(t) \ge~ 0, \linebreak t\ge t_0$.}

See the proof in [13].

In the system (1.1) substitute
$$
\Phi(t) = \exp\biggl\{\int\limits_{t_0}^t\frac{\alpha(\tau)}{p(\tau)} d\tau\biggr\} U(t),
\phantom{aaa} t \ge t_0, \eqno (2.5)
$$
where $\alpha(t)$  is a continuously differentiable function on  $[t_0,+\infty)$, $U(t)\equiv  colon \{u_1(t), u_2(t)\}, \linebreak t\ge~ t_0$.  We obtain
$$
(p(t)U'(t))' + (2\alpha(t) + q(t)) U'(t) + \biggl[\alpha'(t) + \frac{\alpha^2(t)}{p(t)} + \frac{q(t) \alpha(t)}{p(t)}\biggr]U(t) + R(t) U(t) = 0,  \eqno (2.6)
$$
$t\ge t_0$.
Let $\alpha_0(t)$ be a solution of the Riccati equation
$$
\alpha'(t) + \frac{1}{p(t)} \alpha^2(t) + \frac{q(t)}{p(t)} \alpha(t) + r(t) = 0, \phantom{aaa} t\ge t_0, \eqno (2.7)
$$
on $[t_0,+\infty)$.   Then for  $\alpha(t) = \alpha_0(t), \phantom{a} t\ge t_0$,
the system (2.6)  takes the form
$$
(p(t)U'(t))' + [2\alpha_0(t) + q(t)] U'(t) + S(t) U(t) = 0,\phantom{aaa} t\ge t_0, \eqno (2.8)
$$
where  $S(t) \equiv\left(\begin{array}{l} 0 \phantom{aaa} r_1(t)\\
- r_2(t)\phantom{a} 0
\end{array}
\right), \phantom{a}  t\ge t_0$. Set  $\beta_0(t) \equiv \int_{t_0}^t \exp\biggl\{ - \int\limits_{t_0}^\tau\frac{2\alpha_0(s) + q(s)}{p(s)}d s\biggr\}\frac{d\tau}{p(\tau)}, \phantom{a} t\ge t_0$.                      In the system (2.8) substitute
$$
U(t) \equiv V(\beta_0(t)), \phantom{aaa} V(t) \equiv colon\{v_1(t), v_2(t)\}, \phantom{aaa} t\ge t_0. \eqno (2.9)
$$
We come to the system
$$
\beta_0'(\gamma_0(t))^2 V''(t) + S(t) V(t) = 0, \phantom{aaa} t\in [0;\omega), \eqno (2.10)
$$
where  $\gamma_0(t)$ is the inverse function to $\beta_0(t)$,\phantom{a} $\omega\equiv \int\limits_{t_0}^{+\infty} \exp\biggl\{- \int\limits_{t_0}^\tau\frac{2\alpha_0(s) + q(s)}{p(s)}d s\biggr\}\frac{d\tau}{p(\tau)}$.                                 Excluding the unknown  $v_2(t)$  from the last system  we arrive at the scalar equation
$$
\biggl[\frac{\beta_0'(\gamma_0(t))^2}{r_2(\gamma_0(t))}v_1''(t)\biggr]'' + \frac{r_1(\gamma_0(t))}{\beta_0'(\gamma_0(t))^2} v_1(t) = 0,\phantom{aaa} t\in [0;\omega). \eqno (2.11)
$$

{\bf Definition 2.1.} {\it A solution of Eq. (2.7) is called $t_1$-regular if it exists on  $[t_1,+\infty)$.}

{\bf Definition 2.2.}  {\it A $t_1$-regular solution of Eq. (2.7) is called $t_1$-normal if there exists a neighborhood $\stackrel {o} U$ of  $\alpha_0(t)$ such that every solution $\widetilde{\alpha}(t)$  of Eq. (2.7) with  $\widetilde{\alpha}(t_1) \in \stackrel{o}U$ is $t_1$-regular. Otherwise $\alpha(t)$   is called $t_1$-extremal. It is called   lower (upper) $t_1$-extremal solution if every solution $\widetilde{\alpha}(t)$  of Eq. (2.7) with  $\widetilde{\alpha}(t_1) < \alpha(t_1) \phantom{a} (> \alpha_1(t_1))$     is not $t_1$-regular.}

{\bf Definition 2.3.} {\it Eq. (2.7) is called $t_1$-regular if it has a $t_1$-regular solution.}

For brevity, let us introduce some notations that will be needed in the sequel.
$$
I_0\equiv \int\limits_{t_0}^{+\infty}\exp\biggl\{\int\limits_{t_0}^\tau\frac{q(s)}{p(s)}d s \biggr\}r(\tau) d\tau, \phantom{a}    I(t) \equiv \int\limits_t^{+\infty}\exp\biggl\{ - \int\limits_t^\tau\frac{q(s)}{p(s)}d s \biggr\} \frac{d\tau}{p(\tau)},\phantom{aaaaaaaaaaaaaaaaaa}
$$
$$
\phantom{aaaaaaaaaaaaaaaaaaaaaaaaaaaaaa} \nu_u(t)\equiv \int\limits_t^{+\infty}\exp\biggl\{ - \int\limits_t^\tau \frac{2 u(s) + q(s)}{p(s)}d s \biggr\} \frac{d\tau}{p(\tau)}, \phantom{aaa} t\ge t_0.
$$
where $u(t)$ is an arbitrary continuous function on $[t_0,+\infty)$. Note that if for some $t_1 \ge t_0$  the integral   $I(t_1) \phantom{a}(\nu_u(t_1))$ converges then the integral  $I(t) \phantom{a}(\nu_u(t))$ converges too for every $t\ge t_0$.  Denote by  $reg(t_1)$
the set of initial values $\alpha_{(0)}$,  for which every solution $\alpha(t)$ of Eq. (2.7) with  $\alpha(t_0) = \alpha_{(0)}$  is  $t_1$-regular.

{\bf Lemma  2.1.} {\it Assume Eq. (2.7) is   $t_1$-regular. Then  $reg(t_1) = [\alpha_*(t_1),+\infty)$, where  $\alpha_*(t)$ is the unique lower $t_1$-extremal solution of Eq. (2.7).}

See the proof in [13].

In what follows we will assume that Eq. (2.7) is $t_1$-regular for some  $t_1 \ge t_0$. Without loss of generality we may take that $t_1 = t_0$ Then according to Lemma 2.1 Eq. (2.7) has the unique lower $t_0$-extremal solution which in the sequel we will denote always by  $\alpha_*(t)$.

{\bf Remark 2.1.}  {\it The  $t_1$-regularity of Eq.  (2.7) for some $t_1 \ge t_0$ is equivalent to non oscillation of the scalar equation
$$
(p(t)\phi'(t))' + q(t) \phi'(t) + r(t) \phi(t) = 0, \phantom{aaa} t\ge t_0. \eqno (2.12)
$$
(see [14]). Some  $t_0$-regularity criteria for Eq. (2.7) are proved in  [13, 14]. A non oscillation criterion for Eq. (2.12) is proved in  [15]
(e.g. for  $r(t) \le 0, \phantom{a} t\ge t_0$,  Eq. (2.7) is  $t_0$-regular [see [13]]).}

If $\alpha_N(t)$   and $\alpha_{N_1}(t)$ are  $t_0$-normal solutions of Eq. (2.7), then the following relations are valid (see [16]).
$$
\nu_{\alpha_N}(t) < +\infty, \phantom{aaa} t\ge t_0; \eqno (2.13)
$$
$$
\nu_{\alpha_*}(t) = +\infty, \phantom{aaa} t\ge t_0; \eqno (2.14)
$$
$$
\alpha_*(t) = \alpha_N(t) - \frac{1}{\nu_{\alpha_N}(t)}, \phantom{aaa} t\ge t_0; \eqno (2.15)
$$
$$
\int\limits_{t_1}^{+\infty}\frac{|\alpha_N(t) - \alpha_{N_1}(t)|}{p(t)} d t < +\infty, \phantom{aaa}\int\limits_{t_1}^{+\infty}\frac{\alpha_N(t) - \alpha_*(t)}{p(t)} d t = +\infty. \eqno (2.16)
$$

{\bf Lemma  2.2}.{\it  Assume $r(t) \ge 0, \phantom{a} t\ge t_0$, and has an unbounded support. If:

\noindent
$I^*)$ $I(t_0) = +\infty,$  then  $\alpha_*(t) \ge 0, \phantom{a} t\ge t_0$;

\noindent
$II^*)$ $I(t_0) < +\infty,$   then  $\alpha_*(t) <  0, \phantom{a} t\ge t_0$
.}

See the proof in  [11].

{\bf Lemma 2.3.} {\it  Assume  $r(t)$  has an unbounded support and is non negative. Moreover assume   $I_0<~ +\infty$. Then Eq. (2.7) has a positive $t_1$-regular solution for some $t_1 \ge t_0$.}

See the proof in  [14].

{\bf Lemma 2.4}. {\it  Assume  $a(t) \ge 0, \phantom{a} c(t) \le 0, \phantom{a} t\ge t_0$, and have unbounded supports. Then the solution  $y_+(t)$ of Eq. (2.1) with  $y_+(t_0) = 0$  is $t_0$-normal.}

See the proof in [11].

{\bf Lemma 2.5} {\it Assume  $r(t) \le 0, \phantom{a} t\ge t_0$, and has an unbounded support. Then the  $t_0$-extremal solution $\alpha_*(t)$  of Eq. (2.7) is negative.}

Proof. Since by Lemma 2.4 the solution $\alpha_+(t)$  of Eq. (2.7) with   $\alpha_+(t_0)=0$  is $t_0$-normal from Lemma 2.1 it follows that $\alpha_*(t_0) < 0$. Suppose for some $t_1 > t_0 \phantom{a} \alpha_*(t_1) = 0$. Then by virtue of Lemma 2.4  $\alpha_*(t)$ is $t_1$-normal.  As far as the solutions of Eq. (2.7) continuously depend on their initial values then from the relations  $\alpha_*(t_0) < 0, \phantom{a} \alpha_*(t_1) = 0$  and from the  $t_1$-normality of  $\alpha_*(t)$ it is easy to derive that $\alpha_*(t)$ is  $t_0$-normal. The obtained contradiction shows that  $\alpha_*(t) < 0, \phantom{a} t\ge t_0$. The lemma is proved.

{\bf Lemma 2.6.} {\it  Let the following conditions be satisfied:

\noindent
$I^{**}) \phantom{a} r(t) \ge 0; \phantom{a} t \ge t_0$, and has an unbounded support;  $II^{**}) \phantom{a}  I(t_0) < +\infty;$

\noindent
$III^{**})$     Eq. (2.7) has a positive $t_0$-normal solution.

\noindent
Then
$$
\alpha_*(t) > - \frac{1}{I(t)}, \phantom{aaa} t\ge t_0. \eqno (2.17)
$$
}

Proof. Let according to the condition  $III^{**})$  $\alpha_+(t)$  be a positive $t_0$-normal solution of Eq. (2.3). Since
$$
\widetilde{\alpha}_+(t) \equiv \frac{\alpha_+(t_0)\exp\biggl\{ - \int\limits_{t_0}^t\frac{q(s)}{p(s)}ds\biggr\}}{1 + \alpha_+(t_0)\int\limits_{t_0}^t\exp\biggl\{ - \int\limits_{t_0}^\tau\frac{q(s)}{p(s)}ds\biggr\}\frac{d\tau}{p(\tau)}}, \phantom{aaa} t \ge t_0,
$$
is a $t_0$-regular solution of the equation
$$
\alpha'(t) + \frac{1}{p(t)}\alpha^2(t) + \frac{q(t)}{p(t)}\alpha(t) = 0, \phantom{aaa} t\ge t_0,
$$
by Theorem 2.1 from $I^{**})$ it follows that
$$
\alpha_+(t) \le \widetilde{\alpha}_+(t), \phantom{aaa} t\ge t_0. \eqno (2.18)
$$
Note that for every $s\ge t\ge t_0$
$$
\widetilde{\alpha}_+(s) \equiv \frac{\widetilde{\alpha}_+(t)\exp\biggl\{ - \int\limits_t^s\frac{q(\zeta)}{p(\zeta)}d\zeta\biggr\}}{1 + \widetilde{\alpha}_+(t)\int\limits_{t}^s\exp\biggl\{ - \int\limits_{t}^\xi\frac{q(\zeta)}{p(\zeta)}d\zeta\biggr\}\frac{d\xi}{p(\xi)}} \eqno (2.19)
$$
(by the uniqueness theorem the right hand part of (2.19) does not depend on $t$). Since $\alpha_+(t)$ is $t_0$-normal by (2.13)  for every  $t\ge t_0$ the integral $\nu_{\alpha_+}(t)$  converges. From   $II^{**})$,  (2.18) and (2.19) it follows:
$$
\nu_{\alpha_+}(t) \ge \int\limits_t^{+\infty} \exp\biggl\{ - \int\limits_t^\tau\frac{2\widetilde{\alpha}_+(s) + q(s)}{p(s)}d s\biggr\}\frac{d\tau}{p(\tau)} =
$$
$$
=\int\limits_t^{+\infty} \exp\biggl\{ - \int\limits_t^\tau\frac{ q(s)}{p(s)}d s\biggr\}
\exp\biggl\{ - 2 \ln \biggl[1 + \widetilde{\alpha}_+(t) \int\limits_t^\tau \exp\biggl\{- \int\limits_t^\xi \frac{q(\zeta)}{p(\zeta)}d\zeta\biggr\}\frac{d\xi}{p(\xi)}\biggr]\frac{d\tau}{p(\tau)}\biggr\} =
$$
$$
=\frac{1}{\widetilde{\alpha}_+(t)}\int\limits_t^{+\infty}d\left[\frac{1}{1 + \widetilde{\alpha}_+(t) \int\limits_t^\tau \exp\biggl\{- \int\limits_t^\xi \frac{q(\zeta)}{p(\zeta)}d\zeta\biggr\}\frac{d\xi}{p(\xi)}}\right] = \frac{1}{\widetilde{\alpha}_+(t) + \frac{1}{I(t)}}, \phantom{aaa} t\ge t_0. \eqno (2.20)
$$
By  (2.15) the equality $\alpha_*(t) = \alpha_+(t) - \frac{1}{\nu_{\alpha_+}(t)}, \phantom{a} t\ge t_0,$ is satisfied.  This together with (2.18) and (2.20) implies (2.17). The lemma is proved.

{\bf Lemma 2.7}.  {\it Let the following conditions be satisfied:

\noindent
$a) \phantom{a} r(t) \le 0, \phantom{a} t\ge t_0; \phantom{a} b) \phantom{a} I(t_0) < +\infty;$

\noindent
$c) \phantom{a} \int\limits_{t_0}^{+\infty}|r(t)| \biggl(\int\limits_{t_0}^t \exp \biggl\{\int\limits_\tau^t \frac{q(s)}{p(s)} d s \biggr\}\frac{d\tau}{p(\tau)}\biggr) d t < +\infty.$

\noindent
Then the integral $\int\limits_{t_0}^{+\infty}\frac{\alpha_*(\tau)}{p(\tau)}d\tau$ is convergent.}

Proof. Consider the equation
$$
v'(t) - r(t)v^2(t) - \frac{q(t)}{p(t)} v(t) - \frac{1}{p(t)} = 0, \phantom{aaa} t\ge t_0.
$$
Let $v_0(t)$ be a solution of this equation with  $v_0(t_0) > 0$. Then by Lemma 2.1 and Lemma~ 2.4 from a) it follows that  $v_0(t)$ is $t_0$-normal and by virtue of Theorem 2.2
$$
v_0(t_0) > 0, \phantom{aaa} t\ge t_0. \eqno (2.21)
$$
Since  $v_1(t) \equiv v_0(t_0)\exp\biggl\{\int\limits_{t_0}^t\frac{q(s)}{p(s)}d s \biggr\} + \int\limits_{t_0}^t \exp\biggl\{\int\limits_{\tau}^t\frac{q(s)}{p(s)}d s \biggr\}\frac{d\tau}{p(\tau)}, \phantom{a} t\ge t_0,$  is a solution of the equation
$$
v'(t) - \frac{q(t)}{p(t)} v(t) - \frac{1}{p(t)} = 0, \phantom{aaa} t\ge t_0,
$$
In virtue of Theorem 2.1 from a) it follows that
$$
v(t)\le v_0(t_0)\exp\biggl\{\int\limits_{t_0}^t\frac{q(s)}{p(s)}d s \biggr\} + \int\limits_{t_0}^t \exp\biggl\{\int\limits_{\tau}^t\frac{q(s)}{p(s)}d s \biggr\}\frac{d\tau}{p(\tau)}, \phantom{a} t\ge t_0. \eqno (2.21)
$$
Obviously by (2.21) $\alpha_0(t) \equiv \frac{1}{v_0(t)}$  is a $t_0$-regular solution of Eq. (2.7). Then from (2.21) it follows
$$
\frac{1}{\alpha_0(t)} \le \frac{1}{\alpha_0(t_0)}\exp\biggl\{\int\limits_{t_0}^t\frac{q(s)}{p(s)}d s \biggr\} + \int\limits_{t_0}^t \exp\biggl\{\int\limits_{\tau}^t\frac{q(s)}{p(s)}d s \biggr\}\frac{d\tau}{p(\tau)}, \phantom{aaa} t\ge t_0.
$$
Taking into account the condition  a) from here we obtain
$$
\int\limits_{t_0}^t\frac{r(\tau)}{\alpha_0(\tau)} d\tau \ge \frac{1}{\alpha_0(t_0)} \int\limits_{t_0}^t r(\tau)\biggl[\exp\biggl\{\int\limits_{t_0}^\tau \frac{q(s)}{p(s)} d s\biggr\} + \phantom{aaaaaaaaaaaaaaaaaaaaaaaaaaaaaaa}
$$
$$
\phantom{aaaaaaaaaaaaaaaaaaaaaaaaaaaa}+ \alpha_0(t_0)\int\limits_{t_0}^\tau\exp\biggl\{\int\limits_\xi^\tau \frac{q(s)}{p(s)}\biggr\}\frac{d\xi}{p(\xi)}\biggr]d\tau, \phantom{aaa} t\ge t_0. \eqno (2.23)
$$
By the Fubini's theorem from  b) and  c) it follows that the integral  $\int\limits_{t_0}^{+\infty}r(\tau)\exp\biggl\{\int\limits_{t_0}^\tau \frac{q(s)}{p(s)}ds\biggr\}d\tau$  converges. From here from  c)  and from (2.23) it follows
$$
\int\limits_{t_0}^{+\infty}\frac{r(\tau)}{\alpha_0(\tau)}d\tau > - \infty. \eqno (2.24)
$$
Let us prove the equality
$$
\exp\biggl\{\int\limits_{t_0}^t\frac{\alpha_*(s)}{p(s)}d s\biggr\} = \phantom{aaaaaaaaaaaaaaaaaaaaaaaaaaaaaaaaaaaaaaaaaaaaaaaaaaaaaaaaaa}
$$
$$
\phantom{aaa}=\frac{1}{\nu_{\alpha_0}(t_0)}\exp\biggl\{\int\limits_{t_0}^t\frac{\alpha_0(s)}{p(s)}d s\biggr\}\int\limits_t^{+\infty}\frac{\alpha_0(s)}{\alpha_0(t_0)p(s)}
\exp\biggl\{\int\limits_{t_0}^s\biggl[\frac{r(\xi)}{\alpha_0(\xi)} - \frac{\alpha_0(\xi)}{p(\xi)}\biggr]d\xi\biggr\}d s, \eqno (2.25)
$$
$t\ge t_0$.
By (2.15) we have  $\alpha_*(t) = \alpha_0(t) - \frac{1}{\nu_{\alpha_0}(t)}, \phantom{a} t\ge t_0.$                From here it follows
$$
\int\limits_{t_0}^t\frac{\alpha_*(\tau)}{p(\tau)} d\tau = \int\limits_{t_0}^t\frac{\alpha_0(\tau)}{p(\tau)}d\tau - \int\limits_{t_0}^t\frac{1}{p(\tau)\nu_{\alpha_0}(\tau)} d\tau = \int\limits_{t_0}^t\frac{\alpha_0(\tau)}{p(\tau)}d \tau + \phantom{aaaaaaaaaaaaaaaaaaaaaaaaaaaaaaaaaa}
$$
$$
+\int\limits_{t_0}^t d\biggl(\ln\biggl[\int\limits_\tau^{+\infty}\exp\biggl\{- \int\limits_{t_0}^s\biggl[\frac{2 \alpha_0(\xi)}{p(\xi)} + \frac{q(\xi)}{p(\xi)}\biggr]d\xi\biggr\}\frac{d s}{p(s)}\biggr]\biggr) - \ln \nu_{\alpha_0}(t_0), \phantom{aaa} t\ge t_0. \eqno (2.26)
$$
By (2.7) from the inequality $\alpha_0(t) > 0, \phantom{a} t\ge t_0,$ it follows that $\frac{2 \alpha_0(\xi)}{p(\xi)} + \frac{q(\xi)}{p(\xi)} = - \frac{\alpha_0'(\xi) - \frac{\alpha_0^2(\xi)}{p(\xi)} + r(\xi)}{\alpha_0(\xi)}, \linebreak \xi \ge t_0$.  This together with (2.24) implies (2.25). Show that
$$
\int\limits_{t_0}^{+\infty}\frac{\alpha_0(s)}{p(s)}d s = +\infty. \eqno (2.27)
$$
Since
$
\alpha_1(t)\equiv \frac{\alpha_0(t_0)\exp\biggl\{-\int\limits_{t_0}^t\frac{q(s)}{p(s)} ds \biggr\}}{1 + \alpha_0(t_0)\int\limits_{t_0}^t\exp\biggl\{- \int\limits_{t_0}^\tau \frac{q(s)}{p(s)} d s\biggr\}\frac{d\tau}{p(\tau)}}, \phantom{a} t\ge t_0,
$
is a $t_0$-regular solution of the equation
$$
\alpha'(t) + \frac{1}{p(t)}\alpha^2(t) + \frac{q(t)}{p(t)} = 0, \phantom{aaa} t\ge t_0,
$$
by Theorem 2.1 from  a)  it follows that  $\alpha_0(t) \ge \alpha_1(t), \phantom{a} t\ge t_0$.                  Therefore
$\int\limits_{t_0}^t\frac{\alpha_0(\tau)}{p(\tau)} d\tau \ge \linebreak \ge \ln \biggl(1 + \alpha_0(t_0)\int\limits_{t_0}^t\exp\biggl\{- \int\limits_{t_0}^\tau \frac{q(s)}{p(s)} d s\biggr\}\frac{d\tau}{p(\tau)}\biggr), \phantom{a} t\ge t_0
$.
From here and from  b) it follows (2.27). From (2.27) it follows that to the right hand side of  (2.25) we can apply the L'hospitals   rule. Then
$$
\exp\biggl\{\lim\limits_{t\to +\infty}\int\limits_{t_0}^t \frac{\alpha_*(s)}{p(s)} d s\biggr\} = \frac{1}{\nu_{\alpha_0}(t_0)}\lim\limits_{t\to +\infty}\frac{-\frac{\alpha_0(t)}{\alpha_0(t_0)p(t)}
\exp\biggl\{\int\limits_{t_0}^t\biggl[\frac{r(s)}{\alpha_0(s)} - \frac{\alpha_0(s)}{p(s)}\biggr]d s\biggr\}}{- \frac{\alpha_0(t)}{p(t)}\exp\biggr\{-\int\limits_{t_0}^t\frac{\alpha_0(s)}{p(s)} d s\biggr\}} =\phantom{aaaaaaaaaaaaaaaaa}
$$
$
=\frac{1}{\nu_{\alpha_0}(t_0)\alpha_0(t_0)}\lim\limits_{t\to+\infty}\exp\biggl\{\int\limits_{t_0}^t \frac{q(s)}{p(s)} d s\biggr\}.
$
From here and from (2.24) it follows that \linebreak $\exp\biggl\{\lim\limits_{t\to +\infty}\int\limits_{t_0}^t \frac{\alpha_*(s)}{p(s)} d s\biggr\}<  +\infty$. Hence the integral
$\int\limits_{t_0}^{+\infty} \frac{\alpha_*(s)}{p(s)} d s$  converges. The lemma is proved.

\vskip 12pt

{\bf 3. Oscillation criteria}.
 Set:
$$
\beta_*(t) \equiv \int\limits_{t_0}^t\exp\biggl\{ - \int\limits_{t_0}^\tau \frac{2\alpha_*(s) + q(s)}{p(s)} d s\biggr\} \frac{d\tau}{p(\tau)}, \phantom{aaa} t\ge t_0,
$$
and the inverse function of  $\beta_*(t)$ denote by $\gamma_*(t)$.

{\bf Theorem 3.1} {\it  Let the following conditions be satisfied:

\noindent
$A_1) \phantom{a} r(t)  \ge 0, \phantom{a} t\ge t_0, \phantom{a} B_1)\phantom{a} I(t_0) = +\infty;$

\noindent
$C_1) \phantom{a} \int\limits_{t_0}^{+\infty} p(t) |r_k(t)| \biggl(\int\limits_{t_0}^t \exp\biggl\{\int\limits_\tau^t \frac{q(s)}{p(s)} d s\biggr\} \frac{d\tau}{p(\tau)}\biggr) d t = +\infty, \phantom{a} k=1,2.$

\noindent
Then the system $(1.1)$ is oscillatory.}

Proof. By (2.14) we have $\lim\limits_{t\to +\infty} \beta_*(t) = \nu_{\alpha_*}(t_0) = +\infty$. Hence the domain of the function $\gamma_*(t)$ is the half line  $[0,+\infty)$.
From here and from $(2.5)  - (2.11)$  it follows that the theorem will be proved if we show that the equation
$$
\biggl[\frac{\beta_*'(\gamma_*(t))^2}{r_2(\gamma_*(t))}\phi''(t)\biggr]'' + \frac{r_1(\gamma_*(t))}{\beta_*'(\gamma_*(t))^2} \phi(t) = 0, \phantom{aaa} t \ge 0, \eqno (3.1)
$$
is oscillatory. Without loss of generality we may take that $r_k(t)> 0, \ph k=1,2, \ph t\ge t_0$ Then by  Theorem 1.4  the oscillation of (3.1) will be proved if we show that
$$
I_k \equiv \int\limits_0^{+\infty} t \frac{r_k(\gamma_*(t))}{\beta_*'(\gamma_*(t))^2}d t = +\infty, \phantom{aaa} k=1,2. \eqno (3.2)
$$
We have
$$
I_k = \int\limits_{t_1}^{+\infty}\beta_*(t)\exp\biggl\{\int\limits_{t_1}^t \frac{2\alpha_*(s) + q(s)}{p(s)} d s\biggr\} p(t) r_k(t) d t =\phantom{aaaaaaaaaaaaaaaaaaaaaaaaaaaaaaaaaaaaaa}
$$
$$
\phantom{aaaaaaaaaaaa}= \int\limits_{t_1}^{+\infty} p(t) r_k(t) \biggl( \int\limits_{t_1}^t \exp\biggl\{\int\limits_\tau^t \frac{2\alpha_*(s) + q(s)}{p(s)} d s\biggr\}\frac{d\tau}{p(\tau)}\biggr) d t, \phantom{aaa} k=1,2. \eqno (3.3)
$$
By virtue of Lemma 2.2.$I^*$)  from $A_1)$  and $B_1)$  it follows that $\alpha_*(t) \ge 0, \phantom{a} t\ge t_0$.
Then from  (3.3) we obtain
$$
I_k \ge \int\limits_{t_1 + 1}^{+\infty} p(t) r_k(t)  \biggl( \int\limits_{t_1}^t \exp\biggl\{\int\limits_\tau^t\frac{q(s)}{p(s)} d s\biggr\} \frac{d\tau}{p(\tau)}\biggr) d t \stackrel{def}=\widetilde{I}_k, \phantom{a} k=1,2. \eqno (3.4)
$$
Set  $\varepsilon(t) \equiv 1 - \biggl[\int\limits_{t_0}^{t_1}\exp\biggl\{-\int\limits_{t_0}^\tau\frac{q(s)}{p(s)} d s\biggr\} \frac{d\tau}{p(\tau)} \bigg/  \int\limits_{t_0}^{t}\exp\biggl\{-\int\limits_{t_0}^\tau\frac{q(s)}{p(s)} d s\biggr\} \frac{d\tau}{p(\tau)}\biggr], \phantom{a} t\ge t_0.$
We have
$$
\widetilde{I}_k = \int\limits_{t_1 +1}^{+\infty} p(t) r_k(t) \varepsilon(t)  \biggl( \int\limits_{t_0}^t \exp\biggl\{\int\limits_\tau^t \frac{q(s)}{p(s)}d s\biggr\} \frac{d\tau}{p(\tau)}\biggr) d t, \phantom{aaa} k=1,2. \eqno (3.5)
$$
Obviously $\varepsilon(t) \ge \varepsilon(t_0 + 1) > 0$  for $t\ge t_0 + 1$. From here from $C_1)$,  (3.4) and (3.5) it follows (3.2). The theorem is proved.

It is obvious that Theorem 1.2 does not follow from Theorem 3.1. It is also obvious Theorem 3.1 does not follow from Theorem 1.2 as well. Moreover  from (1.5)  is seen that if $p_1(t) = q_2(t), \ph t \ge t_0$, then $\frac{p_2(t)}{q_1(t)} = const$. Therefore in the particular case $p(t) \equiv~ 1, \linebreak q(t)  \equiv~ 0$ Theorem 3.1 does not follow from Theorem 1.2 as well (since the relation $\frac{r_2(t)}{r_1(t)} = const$ may not be satisfied).

{\bf Remark 3.1} {\it Theorem 3.1 is a generalization of Theorem 1.4}

{\bf Remark 3.2}. {\it  The conditions $B_1)$  and  $C_1)$  of Theorem 3.1 are satisfied if in particular  one of the following conditions: $0 < - q(t) \le M = const, \phantom{a} t\ge t_0;$ the function  $\theta(t) \equiv \int\limits_{t_0}^t\frac{q(s)}{p(s)}d s, \phantom{a} t\ge t_0,$ is bounded  is satisfied and the conditions: $\int\limits_{t_0}^{+\infty} \frac{d\tau}{p(\tau)} = \int\limits_{t_0}^{+\infty}p(t) r_k(t) d t =~ +\infty, \linebreak k=1,2,$ are satisfied.}

{\bf Theorem 3.2} {\it Let the condition  $A_1)$ of Theorem 3.1 be satisfied and let

\noindent
$A_2) \phantom{a} r(t)$ has a unbounded support and  $I_0 < +\infty; \phantom{a} B_2) \phantom{a} I(t_0) < +\infty; \linebreak C_2) \int\limits_{t_0}^{+\infty} p(t) |r_k(t)| \exp\biggl\{\int\limits_{t_0}^t\frac{q(s)}{p(s)} d s \biggr\} d t = + \infty, \phantom{a} k=1,2$.
Then the system $(1.1)$ is oscillatory.}

Proof. Without loss of generality we may take that $r_k(t) > 0, \ph k=1,2, \ph t \ge t_0$. Then to prove this theorem it is enough (as in the proof of Theorem 3.1) to show that
$$
I_k \equiv \int\limits_0^{+\infty} t \frac{r_k(\gamma_*(t))}{\beta_*'(\gamma_*(t))^2}d t = +\infty, \phantom{aaa} k=1,2. \eqno (3.6)
$$
Consider the function
$$
J(t) \equiv \exp\biggl\{\int\limits_{t_0}^t\frac{2\alpha_*(s)}{p(s)} d s\biggr\}\int\limits_{t_1}^t \exp\biggl\{- \int\limits_{t_0}^\tau \frac{2\alpha_*(s) + q(s)}{p(s)} d s\biggr\} \frac{d\tau}{p(\tau)}, \phantom{aaa} t \ge t_0.
$$
By Lemma  2.2.$II^*$) from  $B_2)$ it follows that $\alpha_*(t) < 0, \phantom{a} t\ge t_1,$   for some  $t_1 \ge t_0$.   Without loss of generality we can take that $t_1 = t_0$. Then since $\alpha_*(t)$ is negative we have \linebreak
$
J(t) = \exp\biggl\{\int\limits_{t_0}^t \frac{2\alpha_*(s)}{p(s)} d s\biggr\}\times
$
$$
\phantom{aaaaaa}\times \int\limits_{t_0}^t\biggl[\frac{-1}{2\alpha_*(\tau)}\biggr]\biggl(\exp\biggl\{- \int\limits_{t_0}^\tau \frac{2\alpha_*(s)}{p(s)}d s\biggr\}\biggr)' \exp\biggl\{ - \int\limits_{t_0}^\tau\frac{q(s)}{p(s)}d s\biggr\} d \tau, \phantom{aaa} t\ge t_0. \eqno (3.7)
$$
Since $\alpha_*(t) < 0, \phantom{a} t\ge t_0$, we have    $\biggl(\exp\biggl\{- \int\limits_{t_1}^t \frac{2\alpha_*(s)}{p(s)}d s\biggr\}\biggr)'>0, \phantom{a} t\ge t_0.$ By virtue of Lemma~ 2.3  and Lemma  2.6 from here from $A_1), \phantom{a} A_2)$  and  (3.7) it follows that
$$
J(t) \ge \frac{1}{2}\exp\biggl\{\int\limits_{t_0}^t\frac{2\alpha_*(s)}{p(s)}d s\biggr\}\int\limits_{t_0}^t\biggl[\exp\biggl\{- \int\limits_{t_0}^\tau\frac{2\alpha_*(s)}{p(s)} d s\biggr\}\biggr]' I(\tau) \exp\biggl\{- \int\limits_{t_0}^\tau\frac{q(s)}{p(s)} d s\biggr\} d\tau, \phantom{aaa} t\ge t_0.
$$
Using the formula of integration by parts from here we obtain:
$$
J(t) \ge \frac{1}{2}\exp\biggl\{\int\limits_{t_0}^t \frac{2\alpha_*(s)}{p(s)} d s\biggr\}\biggl[\exp\biggl\{- \int\limits_{t_0}^t \frac{2\alpha_*(s)}{p(s)} d s\biggr\} - 1  -\phantom{aaaaaaaaaaaaaaaaaaaaaaaaaaaaa}
$$
$$
\phantom{aaaaaaaaaaaaaaaaaaaaaaaaa}- \int\limits_{t_0}^t\exp\biggl\{- \int\limits_{t_0}^\tau \frac{2\alpha_*(s)}{p(s)} d s\biggr\}\biggl( I(\tau) \exp\biggl\{- \int\limits_{t_0}^\tau\frac{q(s)}{p(s)} d s\biggr\} \biggr)'d\tau\biggr] =
$$
$$
=\frac{1}{2}\biggl[1 - \exp\biggl\{\int\limits_{t_0}^t \frac{2\alpha_*(s)}{p(s)} d s\biggr\}- \int\limits_{t_1}^t\exp\biggl\{- \int\limits_{t_0}^\tau  \frac{2\alpha_*(s)}{p(s)} d s\biggr\}\biggl(\int\limits_\tau^{+\infty}\exp\biggl\{- \int\limits_{t_0}^s \frac{q(\zeta)}{p(\zeta)} d\zeta\biggr\} ds\biggr)'d\tau\biggr] \ge
$$
$
\ge \frac{1}{2}\biggl[1 - \exp\biggl\{\int\limits_{t_0}^t \frac{2\alpha_*(s)}{p(s)} d s\biggr\}\biggr], \phantom{a} t\ge t_0.
$
From here and from the inequality  $\alpha_*(t) < 0, \phantom{a} t\ge t_0$, it follows that
$$
J(t) \ge \frac{1}{2}\biggl[1 - \exp\biggl\{\int\limits_{t_0}^{t_0 + 1} \frac{2\alpha_*(s)}{p(s)} d s\biggr\}\biggr], \phantom{aaa} t\ge t_0+1. \eqno (3.8)
$$
It is easy to show that
$
I_k = \int\limits_{t_0}^{+\infty} p(t) r_k(t) \exp\biggl\{\int\limits_{t_0}^t\frac{q(s)}{p(s)} d s \biggr\}J(t) d t, \phantom{a} k=1,2.
$
From here from $C_2)$  and  (3.8) it follows (3.6). The theorem is proved.

{\bf Remark 3.3.} {\it  The conditions $B_1)$  and $B_2)$ are alternative  each other. Therefore \linebreak Theorem~ 3.1 and Theorem 3.2 complete each other.}

{\bf Theorem 3.3}. {\it  Let for some  $\gamma(t) \in C^1[t_0,+\infty)$ the following conditions be satisfied:

\noindent
$A_3)  \phantom{a}\gamma'(t) + \frac{q(t)}{p(t)} \gamma(t) -\frac{\gamma^2(t)}{p(t)} \le r(t) \le 0, \phantom{a} t\ge t_0$,  where

\noindent
$A_3^0) \phantom{a} \gamma'(t) + \frac{q(t)}{p(t)}\gamma(t) \ge 0, \phantom{a} t\ge t_0, \phantom{a} \gamma(t_0) \ge 0; \phantom{aaa} A_3^1) \phantom{a} \int\limits_{t_0}^{+\infty} \frac{\gamma(t)}{p(t)} d t = +\infty;$

\noindent
$B_3) \phantom{a} \int\limits_{t_0}^{+\infty} p(t) |r_k(t)| \biggl(\int\limits_{t_0}^t \exp\biggl\{\int\limits_\tau^t\frac{q(s) - 4\gamma(s)}{p(s)} d s\biggr\}\frac{d\tau}{p(\tau)}\biggr) d t = +\infty, \phantom{a} k=1,2.
$

\noindent
Then the system $(1.1)$ is oscillatory.}

Proof. Consider the equation
$$
\alpha'(t) + \frac{1}{p(t)}\alpha^2(t) + \frac{q(t)}{p(t)} \alpha(t) + \gamma'(t) +\frac{q(t)}{p(t)} \gamma(t) - \frac{\gamma^2(t)}{p(t)} = 0, \phantom{aaa} t \ge t_0. \eqno (3.9)
$$
Obviously $\alpha_-(t)\equiv - \gamma(t)$ is a $t_0$-regular solution of this equation.
Along with (3.9) consider the equation
$$
\alpha'(t) + \frac{1}{p(t)}\alpha^2(t) + \frac{q(t)}{p(t)} \alpha(t) - \gamma'(t) -\frac{q(t)}{p(t)} \gamma(t) - \frac{\gamma^2(t)}{p(t)} = 0, \phantom{aaa} t \ge t_0. \eqno (3.10)
$$
Obviously $\widetilde{\alpha}_+(t)\equiv \gamma(t)$ is a  $t_0$-regular solution of this equation.
Let  $\alpha_+(t)$ be a solution of Eq. (3.9) with  $\alpha_+(t_0) = \widetilde{\alpha}_+(t_0)$.
Then using Theorem 2.1  to the equations (3.9) and  (3.10) and taking into account  $A_3^0)$  we conclude that  $\alpha_+(t)$ is $t_0$-regular and
$$
\alpha_+(t) \le \gamma(t), \phantom{aaa} t\ge t_0 \eqno (3.11)
$$
By virtue of Theorem 2.2 from $A_3)$ it follows that  $\alpha_+(t) \ge 0, \phantom{a} t\ge t_0$.
From here and from  $A_3^1)$ we obtain  $\int\limits_{t_0}^{+\infty}\frac{\alpha_+(t) - \alpha_-(t)}{p(t)} d t = +\infty$. Then according to (2.16) $\alpha_+(t)$ is $t_0$-normal and $\alpha_-(t)$ is  $t_0$-extremal. Therefore by  (2.15)
$$
\alpha_*(t) = \alpha_+(t) - \frac{1}{\nu_{\alpha_+}(t)}, \phantom{aaa} t \ge t_0.   \eqno (3.12)
$$
Let $\alpha_0(t)$ be a solution of Eq. (2.7) with  $\alpha_0(t_0) = 0$. By  Lemma 2.4   from $A_3)$ $(r(t) \le 0,\phantom{a} t\ge t_0)$ it follows that $\alpha_0(t)$ is $t_0$-normal,
and in virtue of Theorem 2.2
$$
\alpha_0(t) \ge 0, \phantom{aaa} t\ge t_0. \eqno (3.13)
$$
On the strength of Theorem 2.1   from  $A_3)$ it follows that  $\alpha_+(t) \ge \alpha_0(t), \phantom{a} t\ge t_0$.              Then
$$
- \frac{1}{\nu_{\alpha_0}(t)} \ge - \frac{1}{\nu_{\alpha_+}(t)}, \phantom{aaa} t\ge t_0. \eqno (3.14)
$$
Since $\alpha_0(t)$ is $t_0$-normal according to (2.15) we have $\alpha_*(t) = \alpha_0(t) - \frac{1}{\nu_{\alpha_0}(t)}, \phantom{a} t\ge t_0$. From here and from (3.11) - (3.14) it follows:
$$
\alpha_*(t) = \alpha_0(t) - \alpha_+(t) +\alpha_+(t) - \frac{1}{\nu_{\alpha_0}(t)} \ge - \alpha_+(t)  + \alpha_-(t) \ge - 2\gamma(t), \phantom{aaa} t\ge t_0. \eqno (3.15)
$$
Without loss of generality we may take that $r_k(t) > 0, \ph k=1,2, \ph t \ge t_0$.
Then to  prove this theorem (as in the case of the proof of Theorem 3.1) it is enough to prove (3.6). We have:
$
I_k = \int\limits_{t_0}^{+\infty} p(t) r_k(t)  \biggl(\int\limits_{t_0}^t \exp\biggl\{\int\limits_\tau^t\frac{2\alpha_*(s) + q(s)}{p(s)} d s\biggr\}\frac{d\tau}{p(\tau)}\biggr) d t, \phantom{a} k=1,2.
$
From here and from (3.15) it follows:
$
I_k \ge \int\limits_{t_0}^{+\infty} p(t) r_k(t)  \biggl(\int\limits_{t_0}^t \exp\biggl\{\int\limits_\tau^t\frac{ q(s) - 4\gamma(s)}{p(s)} d s\biggr\}\frac{d\tau}{p(\tau)}\biggr) d t, \phantom{a} k=1,2.
$
From here and from  $B_3)$ it follows (3.6). The theorem is proved.

{\bf Remark 3.4}. {\it The conditions  $B_3)$  of Theorem 3.3 are satisfied if in particular \linebreak $\gamma(t) > 0, \phantom{a} \gamma'(t)  \ge 0, \phantom{a} t\ge t_0, \phantom{a} \int\limits_{t_0}^{+\infty}\frac{p(t)r_k(t)}{\gamma(t)}\exp\biggl\{\int\limits_{t_0}^t\frac{q(s)}{p(s)}ds\biggr\}d t = +\infty, \phantom{a} k=1,2,$  the function   $\theta(t) \equiv \int\limits_{t_0}^t \frac{q(s)}{p(s)} d s, \phantom{a} t \ge t_0,$  is bounded from above.}

{\bf Example 3.1}. {\it Assume  $p(t) \equiv 1, \phantom{a} q(t) \equiv 0, \phantom{a} r(t) = - t^2, \phantom{a} 0 < \lambda_0 \le r_k(t) \le \lambda_1, \linebreak k=1,2, \phantom{a} t\ge t_0 =1$. Set:  $\gamma(t) = 2 t, \phantom{a} t\ge 1.$                Then it is easy to check that for this case all conditions of Theorem 3.3 are satisfied whereas the conditions of Theorem 1.4 are not satisfied.}

In the previous oscillation theorems $r(t)$ does not change sign. In contrast of these theorems in the next one $r(t)$ may change sign.

{\bf Theorem 3.4} {\it Let for some  $\gamma(t) \in C^1[t_0,+\infty)$  the following conditions be satisfied:

\noindent
$A_4)  \phantom{a} \gamma'(t) + \frac{1}{p(t)} \gamma^2(t)  + \frac{q(t)}{p(t)}\gamma(t) + r(t) \le 0, \phantom{a} t\ge t_0$;

\noindent
$B_4) \phantom{a} \int\limits_{t_0}^{+\infty}\Bigl|\gamma'(t) + \frac{1}{p(t)}\gamma^2(t) + \frac{q(t)}{p(t)}\gamma(t) + r(t)\Bigr|\biggl(\int\limits_{t_0}^t\exp\biggl\{\int\limits_\tau^t\frac{2\gamma(s) + q(s)}{p(s)} d s\biggr\} \frac{d\tau}{p(\tau)}\biggr) d t < +\infty;$

\noindent
$C_4) \phantom{a} \nu_\gamma(t_0) = +\infty; \phantom{aa} D_4) \phantom{a} \int\limits_{t_0}^{+\infty} p(t) |r_k(t)|   \biggl(\int\limits_{t_0}^t \exp\biggl\{\int\limits_\tau^t \frac{2\gamma(s) + q(s)}{p(s)} d s\biggr\}\frac{d\tau}{p(\tau)}\biggr) d t < +\infty, \phantom{a} k=1,2.$

\noindent
Then the system $(1.1)$ is oscillatory}.

Proof.  As in the proofs of previous theorems here we can take that $r_k(t) > 0, \ph k=~1,2, \linebreak t \ge t_0$. Show that
$$
I_k = \int\limits_0^{+\infty} t \frac{r_k(\gamma_*(t))}{\beta_*'(\gamma_*(t))^2} d t = +\infty, \phantom{aaa} k = 1,2. \eqno (3.16)
$$
From the condition  $A_4)$   it follows that the differential inequality
$$
\eta'(t) + \frac{1}{p(t)}\eta^2(t) + \frac{q(t)}{p(t)}\eta(t) + r(t) \le 0, \phantom{aaa} t\ge t_0,
$$
has a solution on  $[t_0,+\infty)$.  Then (see [14]) Eq. (2.7)  is  $t_0$-regular. By Lemma 2.1 from here it follows that $\alpha_*(t)$ is $t_0$-regular. Therefore
$$
I_k=\int\limits_{t_0}^{+\infty} p(t) r_k(t)\exp\biggl\{\int\limits_{t_0}^t\frac{2\alpha_*(s) + q(s)}{p(s)}d s\biggr\}d t\times \phantom{aaaaaaaaaaaaaaaaaaaaaaaaaaaaaaaaaaaaaaaa}
$$
$$
\phantom{aaaaaaaaaaaaaaaaaaaaaaa}\times\int\limits_{t_0}^t\exp\biggl\{ - \int\limits_{t_0}^\tau\frac{2\alpha_*(s) + q(s)}{p(s)}d s\biggr\}\frac{d\tau}{p(\tau)}, \phantom{aaa} k=1,2.  \eqno (3.17)
$$
In Eq. (2.7) substitute  $\alpha(t) = v(t) + \gamma(t)$. We obtain the equation
$$
v'(t) + \frac{1}{p(t)} v^2(t) + \frac{2\gamma(t) + q(t)}{p(t)} v(t) + r_\gamma(t) = 0, \phantom{aaa} t\ge t_0,
$$
where $r_\gamma(t)\equiv \gamma'(t) + \frac{1}{p(t)}\gamma^2(t) + \frac{q(t)}{p(t)}\gamma(t) + r(t), \phantom{a} t\ge t_0$. By virtue of Lemma 2.1 and Lemma~ 2.4 from  $A_4)$ it follows that this equation has the lower  $t_0$-extremal solution  $v_*(t)$.  Then obviously  $\alpha_*(t) = v_*(t) + \gamma(t), \phantom{a} t\ge t_0$.  From here and from  (3.17) it follows:
$$
I_k=\int\limits_{t_0}^{+\infty} p(t) r_k(t)\exp\biggl\{\int\limits_{t_0}^t\frac{2(v_*(s) + \gamma(s)) + q(s)}{p(s)}d s\biggr\}d t\times \phantom{aaaaaaaaaaaaaaaaaaaaaaaaaaaaaaaaaaaaaaaa}
$$
$$
\phantom{aaaaaaaaaaaaaaaaa}\times\int\limits_{t_0}^t\exp\biggl\{ - \int\limits_{t_0}^\tau\frac{2(v_*(s) + \gamma(s))  + q(s)}{p(s)}d s\biggr\}\frac{d\tau}{p(\tau)}, \phantom{aaa} k=1,2.  \eqno (3.18)
$$
By Lemma 2.5 from  $A_4)$ it follows that $v_*(t) < 0, \phantom{a} t\ge t_0$. From here and from (3.18) we obtain:
$$
I_k\ge\int\limits_{t_0}^{+\infty} p(t) r_k(t)\exp\biggl\{\int\limits_{t_0}^t\frac{2v_*(s)}{p(s)}d s\biggr\} \biggl(\int\limits_{t_0}^t\exp\biggl\{ \int\limits_\tau^t\frac{2\gamma(s) + q(s)}{p(s)}d s\biggr\}\frac{d\tau}{p(\tau)}\biggr) d t,  \eqno (3.19)
$$
$k=1,2.$
By Lemma 2.7 from $A_4)  -  C_4)$ it follows that the integral  $I_3\equiv \int\limits_{t_0}^{+\infty}\frac{2v_*(s)}{p(s)}d s$ is convergent. Then taking into account the inequality  $v_*(t) < 0, \phantom{a} t\ge t_0$, from (3.19) we derive:
$$
I_k\ge\exp\{I_3\}\int\limits_{t_0}^{+\infty} p(t) r_k(t)
\biggl(\int\limits_{t_0}^t\exp\biggl\{ \int\limits_\tau^t\frac{2\gamma(s) + q(s)}{p(s)}d s\biggr\}\frac{d\tau}{p(\tau)}\biggr) d t, \phantom{aaa} k=1,2.
$$
This together with  $D_4)$  implies (3.16). Therefore the system $(1.1)$ is oscillatory. The theorem is proved.

Let us indicate  two particular cases when the conditions of Theorem 3.4 are satisfied.

\noindent
1) $p(t) \in C^1[t_0,+\infty), \phantom{a} -\lambda p'(t) + \lambda^2 p(t) - \lambda q(t)  + r(t) \le 0, \phantom{a}  \lambda = const > 0,$

\noindent
$\phantom{aaaaaaaa}\int\limits_{t_0}^{+\infty}\Bigl|-\lambda p'(t) + \lambda^2 p(t) - \lambda q(t)  + r(t)\Bigr| \biggl(\int\limits_{t_0}^t\exp\biggl\{ - \int\limits_\tau^t\biggl(2\lambda - \frac{q(s)}{p(s)}\biggr)d s\biggr\}\frac{d\tau}{p(\tau)}\biggr) d t < +\infty,$

\noindent
$\int\limits_{t_0}^{+\infty}\exp\biggl\{ - \int\limits_{t_0}^\tau\biggl(2\lambda - \frac{q(s)}{p(s)}\biggr)d s\biggr\}\frac{d\tau}{p(\tau)} = +\infty,
 $

$ \phantom{aaaaaaaaaa}\int\limits_{t_0}^{+\infty} p(t) r_k(t) \biggl( \int\limits_{t_0}^t\exp\biggl\{ - \int\limits_\tau^t\biggl(2\lambda - \frac{q(s)}{p(s)}\biggr)d s\biggr\}\frac{d\tau}{p(\tau)}\biggr) d t = +\infty \phantom{a} (\gamma(t) \equiv - \lambda p(t));$

\noindent
2)\phantom{a}$ q(t) \in C^1[t_0,+\infty), \phantom{a} - q'(t) + r(t) \le 0; \phantom{a} t\ge t_0, $

$\phantom{a}\int\limits_{t_0}^{+\infty} \Bigl|-q(t) + r(t)\Bigr|  \biggl(\int\limits_{t_0}^t \exp\biggl\{- \int\limits_\tau^t\frac{q(s)}{p(s)}d s\biggr\} \frac{d\tau}{p(\tau)}\biggr) d t < +\infty, \phantom{a} \int\limits_{t_0}^{+\infty} \exp\biggl\{\int\limits_{t_0}^\tau\frac{q(s)}{p(s)} d s\biggr\}\frac{d\tau}{p(\tau)} = +\infty,$

\noindent

$ \phantom{aaaaaaaaaaaaaaaaaa}\int\limits_{t_0}^{+\infty} p(t) r_k(t) \biggl(\int\limits_{t_0}^t\exp\biggl\{ - \int\limits_\tau^t\frac{q(s)}{p(s)}d s\biggr\}\frac{d\tau}{p(\tau)}\biggr) d t = +\infty \phantom{a} (\gamma(t) \equiv - q(t)).$

{\bf Example 3.2.} {\it For  $p(t) \equiv 1, \phantom{a} \lambda = 1, \phantom{a} q(t) = 1 + 2 \sin t, \phantom{a} r(t) = 2\sin t - \frac{(\sin t^2)^2}{1 + t^2}, \linebreak t\ge t_0, \phantom{a} \int\limits_{t_0}^{+\infty} r_k(t) d t
= +\infty, \phantom{a} k =1,2,$ the condition 1) is satisfied.}

{\bf Example 3.3.} {\it  For   $p(t) \equiv 1, \phantom{a} q(t) = \sin t, \phantom{a} r(t) = \cos t + \frac{\arctan t^2}{1+ |t^3|}, \phantom{a} t\ge t_0, \phantom{a} \int\limits_{t_0}^{+\infty} t r_k(t) d t =\linebreak = +\infty, \phantom{a} k=1,2,$  the condition 2) is satisfied.}

{\bf Example 3.4.} {\it  For   $p(t) \equiv 1, \phantom{a} q(t) = t, \phantom{a} r(t) = 1 + \frac{\sin e^t}{1 + |t|^\omega}, \phantom{a} \omega > 0, \phantom{a} \int\limits_{t_0}^{+\infty} \frac{r_k(t)}{1 + |t|} d t = +\infty, \linebreak k=1,2,$  the condition 2) is satisfied.}

Note that with the restrictions of the last example the condition  $B_1)$  of Theorem~ 3.1 and the condition $A_3)$
of Theorem 3.2 are not satisfied (although the condition $A_1)$ is fulfilled).
It should be noted also that by the same mentioned above reason of comparison of Theorem 3.1 with Theorem 1.2 we can state that in the particular case $p(t)\equiv 1, \linebreak q(t)\equiv~ 0$ Theorem 3.2, Theorem 3.3 and Theorem 3.4 are not also consequences of Theorem 1.2.

The author is grateful to the Referees whose valuable remarks helped very much to improve the paper.

\vskip 20pt

\centerline{\bf References}

\vskip 20pt

\noindent
1. N. N. Moiseev, Asymptotic Methods in Nonlinear Mechanics. Moscow, ''Nauka'', 1969.

\noindent
2. C. A. Swanson, Comparison and Oscillation Theory of Linear Differential equations.\linebreak \phantom{aa} Academic press. New York and London, 1980.

\noindent
3. W. M. Whyburn, On self - adjoint ordinary differential equations of the fourth order. \linebreak \phantom{aa}  Amer. J Math. 52 (1930), 171 - 196.

\noindent
4. H. Kh. Abdullah, An oscillation criterion for systems of linear second order differential \linebreak \phantom{aa} equations. Journal of Mathematical Sciences: Advances and applications, vol. 5,  \linebreak \phantom{aa} Num. 1, 2010, pp. 93 - 102.

\noindent
5.  W. Leighton and Z. Nehary. On the oscillation of solutions of self - adjoint linear \linebreak \phantom{aaaa}  differential equations of fourth order Trans Amer. Math. Soc. 89(1958), 325 - 388.

\noindent
 6. Sh. Ahmad, On the oscillation of solutions of a class of linear fourth order differential \linebreak \phantom{aaaa}  equations. Pacific Journal of Mathematics, vol. 34, No 2 1970.

\noindent
7. M. S. Kenner. On solutions of Certain Self - Adjoint Differential equations  of Fourth \linebreak \phantom{aaaa} Order. Journal of Mathematical Analysis and Applications, 33, 208 - 305 (1971).

\noindent
8. L Erbe, Hille - Wintner type comparison theorems for self - adjoint fourth   order linear\linebreak \phantom{aaaa} differential equations. Proc. Amer. Math. Soc., vol. 30, Number 3, 1980,  417 - 322.

\noindent
9.  J. Regenda, Oscillation theorems for a class of linear differential equations. Czechos-\linebreak \phantom{aaaa} lovak Mathematical Journal, 34 (109), 1984, pp. 113 - 120.

\noindent
10. O. Polumbiny. On oscillatory solutions of fourth order ordinary differential equations.\linebreak \phantom{aaaa} Czechoslovak Mathematical Journal, 49 (124), 1999, pp. 779 - 790.

\noindent
11. G. A. Grigorian, On the Stability of Systems of Two First - Order Linear Ordinary\linebreak \phantom{aa} Differential Equations, Differ. Uravn., 2015, vol. 51, no. 3, pp. 283 - 292.

\noindent
12.  Grigorian G. A. On Two Comparison Tests for Second-Order Linear  Ordinary\linebreak \phantom{aaa} Differential Equations (Russian) Differ. Uravn. 47 (2011), no. 9, 1225 - 1240; trans-\linebreak \phantom{aaa} lation in Differ. Equ. 47 (2011), no. 9 1237 - 1252, 34C10.

\noindent
13  G. A. Grigorian, "Two Comparison Criteria for Scalar Riccati Equations with\linebreak \phantom{aa} Applications". Russian Mathematics (Iz. VUZ), 56, No. 11, 17 - 30 (2012).

\noindent
14. G. A. Grigorian, Global Solvability of Scalar Riccati Equations. Izv. Vissh.\linebreak \phantom{aa} Uchebn. Zaved. Mat.,vol. 51, 2015, no. 3, pp. 35 - 48.

\noindent
15. G. A. Grigorian, "Some Properties of Solutions to Second - Order Linear Ordinary\linebreak \phantom{aa} Differential Equations". Trudty Inst. Matem. i Mekh. UrO RAN, 19, No. 1, 69 - 80\linebreak \phantom{aa} (2013).

\noindent
16.  G. A. Grigorian, Properties of solutions of Riccati equation, Journal of Contemporary\linebreak \phantom{aa}   Mathematical Analysis, 2007, vol. 42, No 4, pp. 184 - 197.

 \end{document}